\magnification=\magstep1%
\input amstex

\hfuzz=8pt
\NoBlackBoxes
\documentstyle{amsppt}
\hsize=5in
\vsize=7.8in
\vskip5truemm
\nopagenumbers
\hskip 4cm {\it Dedicated to S.\ P.\ Novikov on his 60th birthday.}
\vskip 1cm
\topmatter
\title The Morse-Novikov theory of circle-valued functions
and noncommutative localization \endtitle
\rightheadtext{}
\leftheadtext{}
\abstract{We use noncommutative localization to construct a
chain complex which counts the critical points of a
circle-valued Morse function on a manifold, generalizing the
Novikov complex.  As a consequence we obtain
new topological lower bounds on the minimum number of critical points
of a circle-valued Morse function within a homotopy class, generalizing
the Novikov inequalities.}
\endabstract
\author Michael Farber and Andrew Ranicki\endauthor
\address School of Mathematical Sciences, Tel-Aviv University,
Ramat-Aviv 69978, Israel \endaddress
\address Department of Mathematics
and Statistics, University of Edinburgh, Edinburgh EH9 3JZ, Scotland, UK
\endaddress
\email farber\@math.tau.ac.il, aar\@maths.ed.ac.uk \endemail
\thanks{The research was supported by a grant from the Israel
Academy of Sciences and Humanities, by the Hermann Minkowski Center for
Geometry and by EPSRC grant GR/M20563} \endthanks
\endtopmatter

\define\R{{\Bbb R}}
\define\Z{{\Bbb Z}}
\define\cn{C^{Nov}(M,f)}
\nopagenumbers
\subheading{1. Introduction}
\smallskip
Let us start by recalling the way in which chain complexes are used
to count the critical points of a real-valued Morse function.
\smallskip
A Morse function $f:M \to \R$ on a compact (differentiable)
$m$-dimensional manifold $M$ with $c_i(f)$ critical points of index $i$
determines a handlebody decomposition of $M$ with $c_i(f)$ $i$-handles. The
{\it Morse-Smale complex} is the cellular
chain complex of the corresponding $CW$ decomposition of
the universal cover $\widetilde{M}$ of $M$, a based f.g. free
$\Z[\pi_1(M)]$-module chain complex $C(\widetilde{M})$ with
$$\hbox{\rm rank}_{\Z[\pi_1(M)]}C_i(\widetilde{M})~=~c_i(f)~.$$
For any ring morphism $\rho:\Z[\pi_1(M)] \to R$ there is induced a
based f.g. free $R$-module chain complex
$$C(M;R)~=~R\otimes_{\Z[\pi_1(M)]}C(\widetilde{M})$$
with homology $R$-modules
$H_*(M;R)~=~H_*(C(M;R))$
(which in general depend on $\rho$ as well as $R$).
The number $c_i(f)$
of critical points of index $i$ is bounded from below by
the minimum number $\mu_i(M;R)$ of generators in degree $i$ of a
finite f.g. free $R$-module chain complex which is chain
equivalent to $C(M;R)$
$$c_i(f) \geq \mu_i(M;R)~.$$
For a principal ideal domain $R$ and a ring morphism
$\rho~:~\Z[\pi_1(M)] \to R$
the $R$-coefficient {\it Betti numbers} of $M$ are defined as usual by
$$\eqalign{&b_i(M;R)~=~\hbox{\rm rank}_R\big(H_i(M;R)/T_i(M;R)\big)~,\cr
&q_i(M;R)~=~\text{minimum number of generators of}~T_i(M;R)}$$
with $T_i(M;R) \subseteq H_i(M;R)$ the torsion submodule, and
$$\mu_i(M;R)~=~b_i(M;R) + q_i(M;R) + q_{i-1}(M;R)~.$$
The {\it Morse inequalities}
$$c_i(f) \geq b_i(M;R) + q_i(M;R) + q_{i-1}(M;R)$$
are thus an algebraic consequence of the existence of the Morse-Smale complex.
\smallskip
Now for circle-valued functions.
Given a Morse function $f:M \to S^1$ let $c_i(f)$ denote the
number of critical points of index $i$. The {\it Novikov complex} of \cite{N}
is a based f.g. free chain complex $\cn$ over the principal ideal
domain
$$\Z((z))~=~\Z[[z]][z^{-1}]~,$$
such that
\roster
\item"(i)" $\hbox{\rm rank}_{\Z((z))}C^{Nov}_i(M,f)
=c_i(f)$,
\item"(ii)" $\cn$ is chain equivalent to $C\big(M;\Z((z))\big)$,
with
$$\rho:\Z[\pi_1(M)] @>f_*>> \Z[\pi_1(S^1)]~=~ \Z[z,z^{-1}] \to \Z((z))~.$$
\endroster
The chain complex $\cn$ is constructed geometrically using the gradient flow.
The {\it Novikov inequalities}
$$c_i(f) \geq \mu_i\big(M;\Z((z))\big)~=~
b_i\big(M;\Z((z))\big) + q_i\big(M;\Z((z))\big) + q_{i-1}\big(M;\Z((z))\big)$$
are an algebraic consequence of the existence of a chain complex $\cn$ satisfying
(i) and (ii). The $\Z((z))$-coefficient Betti numbers
are called the {\it Novikov numbers} of $M$. The Novikov numbers
depend only on the cohomology class $\xi=f^*(1) \in H^1(M)$, and so
may be denoted by 
$$b_i\big(M;\Z((z))\big)~=~b_i(\xi)~~,~~q_i\big(M;\Z((z))\big)~=~q_i(\xi)~.$$
\indent
A map $f:M \to S^1$ classifies an infinite cyclic cover
$\overline{M}=f^*\R$ of $M$. We shall assume that $M$ and $\overline{M}$
are connected, so that there is defined a short exact sequence
$$0 \to \pi \to  \pi_1(M) @>f_*>> \pi_1(S^1)=\Z \to 0$$
with $\pi=\pi_1(\overline{M})$.
Let $z \in \pi_1(M)$ be such that $f_*(z)=1$, so that
$$\pi_1(M)~=~\pi\times_{\alpha}\Z~~,~~\Z[\pi_1(M)]~=~\Z[\pi]_{\alpha}[z,z^{-1}]~.$$
with $\alpha:\pi \to \pi;g \mapsto z^{-1}gz$ the monodromy automorphism.
\smallskip
Let $\Sigma$ denote the set of square matrices with entries in $\Z[\pi_1(M)]$ having the form
$1-ze$ where $e$ is a square matrix with entries in $\Z[\pi]$.
A ring morphism $\Z[\pi_1(M)] \to R$ is called {\it $\Sigma$-inverting}
if it sends matrices in $\Sigma$ to invertible matrices over the ring $R$.
There exists 
{\it a noncommutative localization in the sense of  P.\ M.\ Cohn} \cite{C},
a ring $\Sigma^{-1}\Z[\pi_1(M)$ together with a ring morphism
$\Z[\pi_1(M)]\to \Sigma^{-1}\Z[\pi_1(M)]$ which has the universal
property that every $\Sigma$-inverting homomorphism
$\Z[\pi_1(M)] \to R$ has a unique factorization
$$\Z[\pi_1(M)] \to \Sigma^{-1}\Z[\pi_1(M)] \to R~.$$
In particular, the inclusion of the group ring $\Z[\pi_1(M)]$ in {\it the Novikov completion}
$$\widehat{\Z[\pi_1(M)]} ~=~\Z[\pi][[z]]_{\alpha}[z^{-1}]$$
is $\Sigma$-inverting, so that there is a factorization
$$\Z[\pi_1(M)] \to \Sigma^{-1}\Z[\pi_1(M)] \to \widehat{\Z[\pi_1(M)]}~.$$
\indent Pazhitnov \cite{P1} extended the geometric construction of
the Novikov complex to a based f.g. free chain complex $\cn$ over
$\widehat{\Z[\pi_1(M)]}$, such that
\roster
\item"(i)" $\hbox{\rm rank}_{\widehat{\Z[\pi_1(M)]}}C^{Nov}_i(M,f)
=c_i(f)$,
\item"(ii)" $\cn$ is chain equivalent to $C(M;\widehat{\Z[\pi_1(M)]})$,
with
$$\rho~=~\hbox{\rm inclusion}~:~\Z[\pi_1(M)] \to \widehat{\Z[\pi_1(M)]}~.$$
\endroster
Moreover, Pazhitnov \cite{P2},\cite{P3} showed that the geometric construction
of $\cn$ can be adjusted so that the differentials
are rational, in the sense that the entries of their matrices
belong to
$\text{im}(\Sigma^{-1}\Z[\pi_1(M)] \to \widehat{\Z[\pi_1(M)]})$.
If $\pi_1(M)$ is abelian and $\alpha=1$ the localization is a subring
of the completion
$$\Sigma^{-1}\Z[\pi_1(M)]~=~(1+z\Z[\pi])^{-1}\Z[\pi][z,z^{-1}]
\subset \widehat{\Z[\pi_1(M)]}~=~\Z[\pi][[z]][z^{-1}]$$
so that $\cn$ is induced from a chain complex defined
over $\Sigma^{-1}\Z[\pi_1(M)]$. 

Our main result is a direct
construction of such a rational lift of the Novikov
complex, which is valid for arbitrary $\pi_1(M)$ (not necessary abelian) and
arbitrary $\alpha$ :
\smallskip
\proclaim{\bf Main Theorem} For every Morse function $f:M \to S^1$
there exists a based f.g. free $\Sigma^{-1}\Z[\pi_1(M)]$-module chain complex
$\widehat{C}(M,f)$ such that
\roster
\item"(i)" $\hbox{\rm rank}_{\Sigma^{-1}\Z[\pi_1(M)]}\widehat{C}_i(M,f)=c_i(f)$,
\item"(ii)" $\widehat{C}(M,f)$ is chain equivalent to $\Sigma^{-1}C(\widetilde{M})$.
\endroster
\endproclaim
The Main Theorem is proved in section 3 by cutting $M$ at the
inverse image $N=f^{-1}(x) \subset M$ of a regular value $x \in S^1$
of $f$, to obtain a fundamental domain $(M_N;N,zN)$ for $\overline{M}$,
and considering the chain homotopy theoretic properties
of the handle decomposition of the cobordism, with $c_i(f)$ $i$-handles.

The existence of the chain complex $\widehat{C}(M,f)$ has as immediate
consequence:
\proclaim{\bf Corollary (Generalized Novikov inequalities)}
For any Morse function $f:M \to S^1$ and any
$\Sigma$-inverting ring morphism $\rho:\Z[\pi_1(M)] \to R$
$$c_i(f) \geq \mu_i(M;R)~.$$
\endproclaim
Recall that the number $\mu_i(M;R)$ is defined as the minimum number  of generators in
degree $i$ of any
f.g. free $R$-module chain complex, which is chain homotopy equivalent to $C(M;R)$.
The numbers $\mu_i(M;R)$
are homotopy invariants of $M$ as a
space over $S^1$ (via $f$); they depend on the homotopy class of $f$
and on the ring homomorphism $\rho$.
\smallskip
As an example consider the following ring of rational functions
$${\Cal R}~ =~
(1+z\Z[z])^{-1}\Z[z,z^{-1}]$$
introduced in \cite{F}.
This ring is a principal ideal domain \cite{F}. The induced homomorphism
$f_\ast: \pi_1(M) \to \pi_1(S^1)=\Z$ determines a ring homomorphism
$$\rho: \Z[\pi_1(M)]\to {\Cal R}~;~g\mapsto z^{f_\ast(g)}\quad (g\in \pi_1(M))~.$$
It is easy to see that $\rho$ is $\Sigma$-inverting. The corollary implies the inequalities
$$c_i(f)\ge \mu_i(M;\Cal R),\quad i=0, 1, 2, \dots .$$
These inequalities coincide with the classical Novikov inequalities since
$$\mu_i(M;\Cal R)~=~\mu_i\big(M;\Z((z))\big)~=~b_i(\xi)+q_i(\xi)+q_{i-1}(\xi)~.$$
The equivalence between this approach (using the ring $\Cal R$ of rational functions) and the
original approach of Novikov (which used the formal power series ring $\Z((z))$) was proved 
in \cite{F}.

\subheading{2. The endomorphism localization}
\smallskip
This section describes the endomorphism localization of a Laurent polynomial
extension, and provides the algebraic machinery required for the
construction of the chain complex $\widehat{C}(M,f)$ in section 3.
\smallskip
Given a ring $A$ and an automorphism $\alpha:A \to A$
let $A_{\alpha}[z]$, $A_{\alpha}[[z]]$, $A_{\alpha}[z,z^{-1}]$, $A_{\alpha}((z))$ be the
$\alpha$-twisted polynomial extension rings of $A$, with $z$ an
indeterminate over $A$ with $az=z\alpha(a)$ ($a \in A$);
for the record, $A_{\alpha}[z]$ is the ring
of finite polynomials $\sum\limits^{\infty}_{j=0}a_jz^j$, $A_{\alpha}[[z]]$
is the ring of power series $\sum\limits^{\infty}_{j=0}a_jz^j$,
$A_{\alpha}[z,z^{-1}]$ is the ring of finite Laurent polynomials
$\sum\limits^{\infty}_{j=-\infty}a_jz^j$, and $A_{\alpha}((z))=A_{\alpha}[[z]][z^{-1}]$
is the Novikov ring of power series $\sum\limits^{\infty}_{j=-\infty}a_jz^j$
with only a finite number of non-zero coefficients $a_j \in A$ for $j<0$.
\smallskip
\noindent{\bf Definition 2.1} Let $\Sigma$ be the set of
square matrices in $A_{\alpha}[z,z^{-1}]$ of the form $1-ze$ with $e$ a
square matrix in $A$. The {\it endomorphism localization}
$\Sigma^{-1}A_{\alpha}[z,z^{-1}]$ is the noncommutative localization
of $A_{\alpha}[z,z^{-1}]$ inverting $\Sigma$ in the sense of Cohn \cite{C}.\qed
\smallskip
By construction, $\Sigma^{-1}A_{\alpha}[z,z^{-1}]$ is the ring obtained
from $A_{\alpha}[z,z^{-1}]$ by adjoining generators corresponding to the
entries in formal inverses $(1-ze)^{-1}$ of elements $1-ze\in \Sigma$,
and the relations given by the matrix equations
$$(1-ze)^{-1}(1-ze)~=~(1-ze)(1-ze)^{-1}~=~1~.$$
\indent Given an $A$-module $B$ let $z B$ be the $A$-module with elements
$zx$ ($x \in B$) and $A$ acting by
$$A \times z B \to z B ~;~ (a,zx) \mapsto azx~=~z \alpha(a)x~.$$
If $B$ is a f.g. free $A$-module with basis $\{b_1,b_2,\dots,b_r\}$ then
$z B$ is a f.g. free $A$-module with basis $\{zb_1,zb_2,\dots,zb_r\}$.
For f.g. free $A$-modules $B,B'$ an $A_{\alpha}[z,z^{-1}]$-module morphism
$f:B_{\alpha}[z,z^{-1}] \to B'_{\alpha}[z,z^{-1}]$ is a finite
Laurent polynomial
$$f~=~\sum\limits^{\infty}_{j=-\infty}z^jf_j~:~B_{\alpha}[z,z^{-1}] \to B'_{\alpha}[z,z^{-1}]$$
with coefficients $A$-module morphisms $f_j:z^j B \to B'$.
Note that in the special case $\alpha=1:A \to A$ there is defined a natural
$A$-module isomorphism
$$ B \to zB~;~x \mapsto zx~.$$
\proclaim{\bf Proposition 2.2}
(i) A ring morphism $A_{\alpha}[z,z^{-1}] \to R$ which sends every $1-ze \in \Sigma$
to an invertible matrix in $R$ has a unique factorization
$$A_{\alpha}[z,z^{-1}] \to \Sigma^{-1}A_{\alpha}[z,z^{-1}] \to R~.$$
(ii) If $E$ is a f.g. free $A$-module then for any $A$-module morphism
$e:z E \to E$ there is defined
an automorphism of a f.g. free  $\Sigma^{-1}A_{\alpha}[z,z^{-1}]$-module
$$1-ze:\Sigma^{-1}E_{\alpha}[z,z^{-1}] \to \Sigma^{-1}E_{\alpha}[z,z^{-1}]~.$$
(iii) The inclusion $A_{\alpha}[z,z^{-1}] \to A_{\alpha}((z))$ factorizes through $\Sigma^{-1}A_{\alpha}[z,z^{-1}]$
$$A_{\alpha}[z,z^{-1}] \to \Sigma^{-1}A_{\alpha}[z,z^{-1}] \to A_{\alpha}((z))$$
with $A_{\alpha}[z,z^{-1}] \to \Sigma^{-1}A_{\alpha}[z,z^{-1}]$ an injection.
\endproclaim
\demo{Proof}
(i) This is the universal property of noncommutative localization.\newline
(ii) By construction.\newline
(iii) Every matrix of the type $1-ze$ is invertible in $A_{\alpha}((z))$, with
$$(1-ze)^{-1}~=~1+ze+z^2e^2+\dots~.$$
so that (i) applies.\qed
\enddemo
For commutative $A$ and $\alpha=1:A \to A$
$$\Sigma^{-1}A[z,z^{-1}]~=~(1+zA[z])^{-1}A[z,z^{-1}]$$
is just the usual commutative localization inverting the multiplicative subset
$$1+zA[z] \subset A[z,z^{-1}]~.$$
We shall need the following
\proclaim{2.3. Deformation Lemma}
Let $C$ be an $A$-module chain complex of the form
$$d_C~=~\left(\matrix d_D & a & c \cr
0 & d_F & b \cr
0 & 0 & d_{D'} \endmatrix\right)~:~C_i~=~D_i \oplus F_i\oplus D'_i
\to C_{i-1}~=~D_{i-1} \oplus F_{i-1}\oplus D'_{i-1}$$
where $a: F_i\to D_{i-1}$, $b: D'_i\to F_{i-1}$ and $c: D'_i\to D_{i-1}$.
Suppose that the morphism $c: D'_i\to D_{i-1}$ is an isomorphism for all $i$.
The formula
$$\widehat{d}_C~=~d_F-b c^{-1}a~:~F_i \to F_{i-1}$$
defines a "deformed differential" on $F_i$ (i.e. $(\widehat{d}_C)^2=0$),
and the chain complex $\widehat{C}$ defined by
$$\widehat{d}_C~:~\widehat{C}_i~=~F_i \to \widehat{C}_{i-1}~=~F_{i-1}$$
is chain equivalent to $C$.
\endproclaim
\demo{Proof} Note that ${(d_C)}^2=0$ implies
$${(d_D)}^2~=~0~,~{(d_F)}^2~=~0~,~{(d_{D'})}^2~=~0~,$$
and also
$$
\aligned
&d_Da +a d_F~=~0~,\\
&d_Fb +b d_{D'}~=~0~,\\
&d_Dc +c d_{D'}+ab ~=~0~.
\endaligned\tag1
$$
Using (1) we obtain
$$\aligned
(\widehat{d}_C)^2~
&=~(d_F -bc^{-1}a)\cdot (d_F -bc^{-1}a)\\
&=~ - d_Fbc^{-1}a - bc^{-1}a d_F +bc^{-1}a\cdot bc^{-1}a\\
&=~- d_Fbc^{-1}a - bc^{-1}a d_F -
bc^{-1}[d_Dc +c d_{D'}]c^{-1}a\\
&=~0~.
\endaligned$$
Now we define two chain maps:
$$u~=~\left(\matrix -bc^{-1} & 1 & 0  \endmatrix\right): C \to \widehat{C} ,\quad
\text{and}\quad v~=~\left(\matrix 0\\ 1\\ 
-c^{-1}a\endmatrix\right): \widehat{C} \to C~.$$
One checks that
$$ud_C = \widehat{d}_Cu ,\quad d_C v = v \widehat{d}_C, \quad uv = 1_{\widehat{C}}$$
and
$$vu = 1_C - d_C w -wd_C,\quad\text{where}\quad w\, =\, \left(\matrix 0 & 0 & 0\\
0 & 0 & 0\\
c^{-1} & 0 & 0\endmatrix\right): C \to C.
$$
Hence $u$ and $v$ are mutually inverse chain equivalences.\qed
\enddemo
\proclaim{\bf Theorem 2.4}
Let $C$ be a finite based f.g. free $A_{\alpha}[z,z^{-1}]$-module chain complex of the form
$$C~=~{\Cal C}(g-zh:D_{\alpha}[z,z^{-1}] \to E_{\alpha}[z,z^{-1}])$$
where $g:D \to E$, $h:z D \to E$ are
chain maps of finite based f.g. free $A$-module chain complexes, and
${\Cal C}$ denotes the algebraic mapping cone.
If each $g:D_i \to E_i$ is a split injection sending basis elements to
basis elements, then there is defined a $\Sigma^{-1}A_{\alpha}[z,z^{-1}]$-module chain complex
$\widehat{C}$ such that
\roster
\item"(i)" each $\widehat{C}_i$ is a based f.g. free
$\Sigma^{-1}A_{\alpha}[z,z^{-1}]$-module with
$$\hbox{\rm rank}_{\Sigma^{-1}A_{\alpha}[z,z^{-1}]}\widehat{C}_i~=~
\hbox{\rm rank}_AE_i - \hbox{\rm rank}_AD_i~,$$
\item"(ii)" there is defined a chain equivalence
$\Sigma^{-1}C \to \widehat{C}$.
\endroster
\endproclaim
\demo{Proof} Let $F_i \subseteq E_i$ be the submodule generated
by the basis elements not coming from $D_i$, so that
$$\aligned
&g~=~\left(\matrix 1 \cr 0 \endmatrix\right)~:~
D_i \to E_i~=~D_i \oplus F_i~,\cr
&h~=~\left(\matrix e \cr f \endmatrix\right)~:~
z D_i \to E_i~=~D_i \oplus F_i~,\cr
&d_E~=~\left(\matrix d_D & a \cr 0 & d_F \endmatrix\right)~:~
E_i~=~D_i \oplus F_i \to E_{i-1}~=~D_{i-1} \oplus F_{i-1}~,\cr
&d_C~=~\left(\matrix d_D & a & 1-ze \cr
0 & d_F & -zf \cr
0 & 0 & d_D \endmatrix\right)~:\cr
&\hskip20pt
C_i~=~(D_i \oplus F_i\oplus D_{i-1})_{\alpha}[z,z^{-1}]
\to C_{i-1}~=~(D_{i-1} \oplus F_{i-1}\oplus D_{i-2})_{\alpha}[z,z^{-1}]~.
\endaligned
$$
Now apply Lemma 2.3 to the induced chain complex $\Sigma^{-1}C$ over
$\Sigma^{-1}A_{\alpha}[z,z^{-1}]$ with
$$\aligned
&d_{\Sigma^{-1}C}~=~\left(\matrix d_D & a & 1-ze \cr
0 & d_F & -zf \cr
0 & 0 & d_D \endmatrix\right)~:\cr
&\hskip20pt
\Sigma^{-1}C_i~=~\Sigma^{-1}(D_i \oplus F_i\oplus D_{i-1})_{\alpha}[z,z^{-1}]  \cr
&\hskip100pt \to \Sigma^{-1}C_{i-1}~=~\Sigma^{-1}(D_{i-1} \oplus F_{i-1}\oplus D_{i-2})_{\alpha}[z,z^{-1}]
\endaligned
$$
where each
$$c~=~1-ze~:~\Sigma^{-1}(D_{i-1})_{\alpha}[z,z^{-1}] \to
\Sigma^{-1}(D_{i-1})_{\alpha}[z,z^{-1}]$$
is an automorphism. Explicitly, Lemma 2.3 gives a based f.g. free
$\Sigma^{-1}A_{\alpha}[z,z^{-1}]$-module chain complex $\widehat{C}$ with
$$\widehat{d}_C~=~d_F+(zf)(1-ze)^{-1}a~:~
\widehat{C}_i~=~\Sigma^{-1}(F_i)_{\alpha}[z,z^{-1}] \to
\widehat{C}_{i-1}~=~\Sigma^{-1}(F_{i-1})_{\alpha}[z,z^{-1}]$$
which satisfies (i) and (ii).\qed
\enddemo
\subheading{3. The chain complex $\widehat{C}(M,f)$}
\smallskip
Given a Morse function $f:M \to S^1$ we now construct a
chain complex $\widehat{C}(M,f)$ over $\Sigma^{-1}\Z[\pi_1(M)]$
satisfying the conditions of the Main Theorem.
\smallskip
Choose a regular value $x \in S^1$ for $f:M\to S^1$, and cut $M$ along
the codimension 1 framed submanifold
$$N^{m-1}~=~f^{-1}(x) \subset M^m$$
to obtain a fundamental domain $(M_N;N,zN)$ for the infinite cyclic cover
$$\overline{M}~=~\bigcup\limits^{\infty}_{j=-\infty}z^jM_N$$
of $M$, with a Morse function
$$f_N~:~(M_N;N,zN) \to ([0,1];\{0\},\{1\})$$
such that $f_N$ has exactly as many critical points of index $i$ as $f$
$$c_i(f_N)~=~c_i(f)~.$$
Let $\widetilde{N}$, $\widetilde{M}_N$ be the covers of $N$, $M_N$ obtained from
the universal cover $\widetilde{M}$ of $M$ by pullback along the inclusions
$N \to \overline{M}$, $M_N \to \overline{M}$.
The cobordism $(M_N;N,zN)$ has a handle decomposition
$$M_N~=~N \times I \cup \bigcup^m_{i=0}\bigcup\limits_{c_i(f)}h^i$$
with $c_i(f)$ $i$-handles $h^i=D^i \times D^{m-i}$.
Let $\pi_1(\overline{M})=\pi$ with monodromy automorphism $\alpha:\pi \to \pi$,
so that $\pi_1(M)=\pi \times_{\alpha}\Z$, $\Z[\pi_1(M)]=\Z[\pi]_{\alpha}[z,z^{-1}]$
as in the Introduction. The relative cellular
chain complex $C(\widetilde{M}_N,\widetilde{N})$ is a based f.g. free
$\Z[\pi]$-module chain complex with
$$\hbox{\rm rank}_{\Z[\pi]}C_i(\widetilde{M}_N,\widetilde{N})~=~c_i(f)~.$$
Choose an arbitrary $CW$ structure for $N$, and let $c_i(N)$ be the number
of $i$-cells. Let $M_N$ have the $CW$ structure with $c_i(N)+c_i(f)$
$i$-cells:
for each $i$-cell $e^i \subset N$ there is defined an $i$-cell $e^i \times I \subset M_N$,
and for each $i$-handle $h^i$ there is defined an $i$-cell $h^i \subset M_N$.
Then $M=M_N/(N=zN)$ has a $CW$ complex structure with
$c_i(f) + c_i(N) + c_{i-1}(N)$ $i$-cells.
The universal cover $\widetilde{M}$ has cellular
$\Z[\pi_1(M)]$-module chain complex
$$C(\widetilde{M})~=~{\Cal C}(g-zh:C(\widetilde{N})_{\alpha}[z,z^{-1}] \to C(\widetilde{M}_N)_{\alpha}[z,z^{-1}])~
$$
where 
$$g~:~C(\widetilde{N}) \to C(\widetilde{M}_N)~~,~~
h~:~z C(\widetilde{N}) \to C(\widetilde{M}_N)$$
are the $\Z[\pi]$-module chain maps induced by the inclusions
$$g~:~N \to M_N~~,~~h~:~zN\to M_N~.$$
Since $g:N\to M_N$ is the inclusion of a subcomplex
$$g~=~\left(\matrix 1 \cr 0\endmatrix\right)~:~
C_i(\widetilde{N}) \to C_i(\widetilde{M}_N)~=~
C_i(\widetilde{N})\oplus C_i(\widetilde{M}_N,\widetilde{N})$$
is a split injection.  Now apply Theorem 2.4 to the based f.g. free
$\Z[\pi_1(M)]$-module chain complex
$$C(\widetilde{M})~=~{\Cal C}(
g-zh:C(\widetilde{N})_{\alpha}[z,z^{-1}] \to C(\widetilde{M}_N)_{\alpha}[z,z^{-1}])$$
with $D=C(\widetilde{N})$, $E=C(\widetilde{M}_N)$.\qed

\enddocument